\documentclass{article} 
\newtheorem{theorem}{Theorem}

\newtheorem{corollary}[theorem]{Corollary}

\newtheorem{definition}[theorem]{Definition}
\newtheorem{example}[theorem]{Example}

\newtheorem{lemma}[theorem]{Lemma}

\newtheorem{proposition}[theorem]{Proposition}
\newtheorem{remark}[theorem]{Remark}

\newenvironment{proof}[1][Proof]{\noindent\textbf{#1.} }{\ \rule{0.5em}{0.5em}}
\usepackage{amssymb}
\usepackage{amsmath}
\usepackage{indentfirst}
\usepackage[all]{xy}
\usepackage{hyperref}

\numberwithin{subsection}{section}
\numberwithin{equation}{section}

\begin{document}
\title{On Coalgebras over Algebras}
\author{Adriana Balan\thanks{Supported by a Royal Society International Travel Grant} \\
Department of Mathematics II\\ University Politehnica of Bucharest \and Alexander Kurz \\
Department of Computer Science\\ University of Leicester}
\maketitle
 
\begin{abstract} 
We extend Barr's well-known characterization of the final coalgebra of a $Set$-endofunctor as the completion of its initial algebra to the Eilenberg-Moore category of algebras for a $Set$-monad $\mathbf{M}$ for functors arising as liftings. As an application we introduce the notion of commuting pair of endofunctors with respect to the monad $\mathbf{M}$ and show that under reasonable assumptions, the final coalgebra of one of the endofunctors involved can be obtained as the free algebra generated by the initial algebra of the other endofunctor.
\end{abstract}

\section{Introduction}

For any category $\mathcal{C}$ and any $\mathcal{C}$-endofunctor $H$,
there is a \linebreak canonical arrow between the least and the
greatest fixed points of $H$, namely between its initial algebra and
final coalgebra, assuming these exist. Functors for which these
objects exist and coincide were called algebraically compact by Barr
\cite{barr1} - for example, if the base category is enriched over
complete metric spaces \cite{ar} or complete partial orders \cite{sp},
then mild conditions ensure that the endofunctors are algebraically
compact.  However, if the category lacks any enrichment, as $Set$,
this coincidence does not happen. But there is still something to be
said: Barr \cite{barr} showed that for bicontinuous
$Set$-endofunctors, the final coalgebra can be realized as the
completion of its initial algebra. But this works if the functor does
not map the empty set into itself, otherwise the initial algebra would
be empty. Hence some well-known examples are lost, like functors
obtained from powers and products. Barr's result was extended to all
locally finitely presentable categories by Ad{\'a}mek \cite{adamek
  ideal,adamek final}, in the sense that the completion procedure
works for hom-sets, not objects, with respect to all finitely
presentable objects.

In the present paper we have focused on coalgebras the carriers of
which are algebras for a $Set$-monad, not necessarily finitary (see
for example \cite{corr-etal:cmcs98}, \cite{turi-plot:sos}). Our
interest arises from the following two developments.

First, streams or weighted automata, as pioneered by Rutten (\cite%
{rutten:automata}, \cite{rutten:diff-equ}, \cite{rutten:counting}) are
mathematically highly interesting examples of coalgebras, despite the
fact that the type functor is very simple, just $HX=A\times X$ in the
case of streams. The interesting structure arises from $A$, which in
typical examples carries the structure of a semi-ring. In this paper,
we shall bring this structure to the fore by lifting $H$ to the
category of modules for a semi-ring, or more generally, to the
category of algebras for a suitable monad.

Second, in recent work of Kissig and the second author \cite{trace},
it turned out that it is of interest to move the trace-semantics of
Hasuo-Jacobs-Sokolova \cite{HJS:cmcs06} from the Kleisli-category of
a commutative monad to the Eilenberg-Moore category of algebras (for
example, this allows to consider wider classes of monads). Again, for
trace semantics, semi-ring monads are of special interest.

In the first part of this paper, we show that Barr's theorem
\cite{barr} extends from coalgebras on $Set$ to coalgebras on the
Eilenberg-Moore category of algebras $Alg(\mathbf{M})$ for a monad
$\mathbf{M}$ on $Set$, dropping the assumption $H0\not=0$ (hence
allowing examples like the functor $H$ of stream coalgebras mentioned
above).

We consider the situation of a $Set$-endofunctor $H$ that has a
lifting to $Alg(\mathbf{M})$. Under some reasonable assumptions, we
are able to prove that the final $H$-coalgebra can be obtained as the
Cauchy completion of the image of the initial algebra for the lifted
functor, with respect to the usual ultrametric inherited from the
final sequence. For this, we need to understand better the initial
algebra of the lifted functor. This is the purpose of the second part
of the paper, where the special case of an initial algebra which is
free (as an $\mathbf{M}$-algebra) is exhibited. Namely, for two
endofunctors $H$, $T$ and a monad $\mathbf{M}$ on $Set$, we call
$(T,H)$ an $\mathbf{M}$-commuting pair if there is a natural
isomorphism $HM\cong MT$, where $M$ is the functor part of the
monad. It follows that if both algebra lift of $H$ and Kleisli lift of $T$ exist, then
mild requirements ensure that $\widetilde{H}$, the algebra lifted functor of $H$, is
equivalent with the extension of $T$ to
$Alg(\mathbf{M})$ if and only if they form a commuting pair. If this
is the case, then one can recover the initial algebra for the lifted
endofunctor $\widetilde{H}$ as the free $\mathbf{M}$-algebra built on
the initial $T$-algebra.

\section{Final coalgebra for endofunctors lifted to categories of algebras}

\subsection{Final sequence for Set-endofunctors\label{terminal sequence}}

Consider an endofunctor $H:Set\longrightarrow Set$. From the unique arrow \linebreak $%
t:H1\longrightarrow 1$ we may form the sequence

\begin{equation}\label{terminal sequence}
  \xymatrix{
    1 & H1 \ar[l]_t & \ldots \ar[l] & H^n1 \ar[l] & H^{n+1}1 \ar[l]_{H^nt} & \ldots \ar[l] }
\end{equation}
\noindent
Denote by $L$ its limit, with $p_{n}:L\longrightarrow H^{n}1$ the
corresponding cone. As we work in $Set$, recall that the limit $L$ can be
identified with a subset of the cartesian product $\prod\limits_{n\geq 0}H^{n}1$, namely%
\begin{equation*}
L=\{(x_{n})_{n\geq 0}\mid H^{n}t(x_{n+1})=x_{n}\}
\end{equation*}%

By applying $H$ to the sequence and to the limit, we get a cone 
\begin{equation*}
  \xymatrix{1 & H1 \ar[l]_t & \ldots \ar[l] & H^n1 \ar[l] & \ldots \ar[l]_{H^nt} \\
	& & & L \ar@/^1pc/[ulll] \ar[u]^{p_n} & & \\
	& & & HL \ar@{.>}[u]^{\tau} \ar@/^1pc/[uulll] \ar@/_1.3pc/[uu]_{Hp_{n-1}} & &  
}
\end{equation*}
\noindent
with $HL\rightarrow 1$ the unique map to the singleton set. The limit property leads to a map  $\tau :HL\rightarrow L$ such that $p_n \circ \tau=Hp_{n-1}$.

For each $H$-coalgebra $(C,\xi _{C}:C \longrightarrow HC)$ it exists a cone $\alpha
_{n}:C\longrightarrow H^{n}1$ over the sequence (\ref{terminal sequence}),
built inductively as follows: $\alpha _{0}:C\longrightarrow 1$ is the unique
map, then if $\alpha _{n}:C\longrightarrow H^{n}1$ is already obtained,

construct $\alpha _{n+1}$ as the composite
\begin{equation}
C\overset{\xi _{C}}{\longrightarrow }HC\overset{H\alpha _{n}}{%
\longrightarrow }H^{n+1}1  \label{cone C-->Hn1}
\end{equation}%
Then the unique map $\alpha _{C}:C\longrightarrow L$ such that $p_{n} \circ \alpha _{C}=\alpha _{n}$ satisfies the following algebra-coalgebra diagram (\cite{moss}):
\begin{equation*}
  \xymatrix{
    C \ar[r]^{{\alpha}_C} \ar[d]_{{\xi}_C} & L  \\
    HC \ar[r]^{H{\alpha}_C} & HL \ar[u] _{\tau}
}
\end{equation*}

On the sequence (\ref{terminal sequence}), endow each set $H^{n}1$ with the
discrete topology (so all maps $H^{n}t$ will be continuous). Then put the
initial topology \cite{schaefer} coming from this sequence on $L$ and $HL$. It follows that $%
\tau $ is continuous. In particular, the topology on $L$ is given by an
ultrametric: the distance between any two points in $L$ is $2^{-n}$, where $%
n $ is the smallest natural number such that $p_{n}(x)\neq p_{n}(y)$. The
cone $\alpha _{n}:C\longrightarrow H^{n}1$ yields on any coalgebra a
pseudo-ultrametric (hence a topology) and the unique map $%
\alpha _{C}:C\longrightarrow L$ is continuous with respect to it.

If $H$ is $\omega ^{op}$-continuous, it preserves the limit $L$, hence
the isomorphism $\xi={{\tau}^{-1}} :L\simeq HL$ makes $L$ the final $H$-coalgebra. Moreover, using the above topology, the map $\xi$ is a homeomorphism and verifies
\begin{equation}\label{projections from the limit}
Hp_{n-1} \circ {\xi}=p_n
\end{equation}

\subsection{Lifting to Eilenberg-Moore category of algebras for a monad\label%
{lifting}}

Let $\mathbf{M}=(M,M^{2}\overset{m}{\longrightarrow } M,Id\overset{u}{\longrightarrow } M)$ be a a
monad on $Set$. Denote by $Alg({\mathbf{M}})$ the Eilenberg-Moore category of 
$\mathbf{M}$-algebras and by $F^{\mathbf{M}}\dashv U^{\mathbf{M}}:Alg({\mathbf{M}})\longrightarrow Set$ the adjunction between the free and the
forgetful functor. Then $Alg({\mathbf{M}})$ has an initial object, namely $(M0,M^20\overset{m_0}{\longrightarrow } M0)$, the free algebra on the empty set, and a terminal object $1$, the singleton, with algebra structure given by the unique map 
$M1\longrightarrow1$.

For a $Set$-endofunctor  $H$, it is well known (\cite{johnstone}) that
liftings of $H$ to $Alg({\mathbf{M}})$, i.e. endofunctors $\widetilde{H}$ on $%
Alg({\mathbf{M}})$ such that the diagram%
\begin{equation}\label{algebra lift}
\xymatrix{
	Alg({\mathbf{M}}) \ar[r]^{\widetilde{H}} \ar[d]_{U^{\mathbf{M}}} & 
	Alg({\mathbf{M}}) \ar[d]_{U^{\mathbf{M}}} \\
	Set \ar[r]^H & Set
}
\end{equation}
\noindent
commutes, are in one-to-one correspondence with natural transformations \linebreak%
$\lambda :MH\longrightarrow HM$
satisfying
\begin{equation}\label{distributive law endofunctor monad}
	\xymatrix{ H \ar[r]^{u_H} \ar[dr]_{Hu} & MH \ar[d]^{\lambda} \\
	 & HM}
\qquad
	\xymatrix {M^{2}H \ar[r]^{M\lambda} \ar[d]_{m_H} & MHM \ar[r]^{\lambda _M} & 
	HM^2 \ar[d]^{Hm} \\
	MH \ar[rr]^{\lambda} & & HM}
\end{equation} 
\begin{remark} It is worth noticing that the lifting is not unique (as
  there may be more than one distributive law $\lambda:MH
  \longrightarrow HM$). For example, take $G$ a group and
  $HX=MX=G\times X$; consider $H$ as an endofunctor and $M$ as a monad
  with natural transformations $u,m$ obtained from the group
  structure. The algebras for this monad are the $G$-sets. Then it is
  easy to see that a map $f:G \times G \longrightarrow G \times G$
  induces a distributive law $\lambda:MH \longrightarrow HM$ if it
  satisfies $f(e,x)=(x,e)$ for all $x \in G$, where $e$ stands for the
  unit of the group, and $f(\mu \times G)=(G \times \mu)(f \times G)(G
  \times f)$, where we have denoted by $\mu$ the group
  multiplication. Take now $f_1(x,y)=(xy,x)$ and
  $f_2(x,y)=(xyx^{-1},x)$; these maps produce two distributive laws
  $\lambda_1,\lambda_2:MH \longrightarrow HM$ which do not give same lifting
  $\widetilde{H}$, as the $G$-action on $HX$ would be $(x,y,z)
  \longrightarrow (xy,x\rightharpoonup z)$ for $\lambda_1$,
  respectively $(x,y,z) \longrightarrow (xyx^{-1},x\rightharpoonup z)$
  for $\lambda_2$. Here $x,y \in G, z \in X$ and $\rightharpoonup$
  denotes the left $G$-action on $X$. If the liftings would be isomorphic, then the associated categories of coalgebras should also
  be isomorphic. In particular, notice that $H$ is a comonad (as any
  set, in particular $G$, carries a natural comonoid structure) and
  both maps $f_1,f_2$ are actually inducing monad-comonad distributive laws
  $\lambda_1$, respectively $\lambda_2$. Hence each lifting \linebreak carries a comonad structure such that the associated categories of coalgebras for the lifted functors are Eilenberg-Moore categories of coalgebras and they should also be isomorphic. But for $f_1$, a corresponding coalgebra is the same as a $G$-set $(X,\rightharpoonup)$ endowed with
  a map $\theta:X\longrightarrow G$ such that $\theta(g
  \rightharpoonup x)=g\theta(x)$, while for the second structure, the
  compatibility relation yields a crossed $G$-set, i.e. $\theta(g
  \rightharpoonup x)=g\theta(x)g^{-1}$.
\end{remark}
Assume from now on that a lifting of $H$ to $Alg(\mathbf{M})$ exists, given by $\lambda :MH\longrightarrow HM$. For any $\mathbf{M}$-algebra $(X,x)$, $HX$ becomes an
algebra with \linebreak $\xymatrix{MHX \ar[r]^{\lambda _X} & HMX \ar[r]^{Hx} & HX}$ and for any algebra map $(X,x)\longrightarrow (Y,y)$, the corresponding arrow $HX\longrightarrow HY$ respects the algebra structure. Also, for any $H$-coalgebra $(C,C\overset{\xi _{C}}{\longrightarrow }HC)$%
, $MC$ inherits an $H$-coalgebra structure by $\xymatrix{\xi :MC \ar[r]^{M\xi _{C}} & MHC \ar[r]^{\lambda _C} & HMC}$. In particular, if the final coalgebra $(L,L\overset{\xi }{\longrightarrow }HL)$ exists, then
there is a unique coalgebra map $\gamma :ML\longrightarrow L$, given by:
\begin{equation}\label{ML=H-coalgebra}
  \xymatrix{ ML \ar@{.>}[d]_{\gamma} \ar[r]^{M\xi} & MHL \ar[r]^{\lambda _L} & HML \ar@{.>}[d]^{H \gamma}\\
L \ar[rr]^{\xi} & & HL
  	}
\end{equation}

\noindent
Then $(L,\gamma )$ and $(HL,H\gamma \lambda _{L})$ are $\mathbf{M}$-algebras
and $\xi :(L,\gamma )\longrightarrow (HL,H\gamma \lambda _{L})$ becomes an $%
\mathbf{M}$-algebra map. By the lifting property, $\widetilde{H}(L,\gamma
)=(HL,H\gamma \lambda _{L})$ and as any $\widetilde{H}$-coalgebra (its
underlying set) is the carrier of an $H$-coalgebra, it follows that $%
((L,\gamma ),\xi )$ is the final $\widetilde{H}$-coalgebra. Hence despite the fact that the lifting might not be unique, the underlying set of the final $H$-coalgebra is preserved (but with possibly different algebra structure, depending on $\lambda$).

Coming back to the final sequence (\ref{terminal sequence}), note that
any term $H^{n}1$ is an $\mathbf{M}$-algebra by:

\begin{itemize}
\item the obvious unique $\mathbf{M}$-algebra structure on $1$, $%
a_{0}:M1\longrightarrow 1$;

\item given $a_{n}:MH^{n}1\longrightarrow H^{n}1$, define $a_{n+1}$ as the composite
\begin{equation}
MH^{n+1}1\overset{\lambda _{H^{n}1}}{\longrightarrow }HMH^{n}1%
\overset{Ha_{n}}{\longrightarrow }H^{n+1}1  \label{Hn1=M-algebra}
\end{equation}
\end{itemize}
Moreover, all maps in the sequence (\ref{terminal sequence}) are
$\mathbf{M}$%
-algebra maps by (\ref{distributive law endofunctor monad}). \linebreak Applying
$M$ to the sequence produces a cone from $ML$. If we assume $H$ $\omega^{op}$-continuous (hence $\xi :L \simeq HL$ is an isomorphism), we can understand better this cone-construction:

\begin{lemma}
The cone $(ML\overset{Mp_{n}}{\longrightarrow }MH^{n}1\overset{a_{n}}{%
\longrightarrow }H^{n}1)_{n\geq 0}$ coincides with the cone $\alpha
_{n}:ML\longrightarrow H^{n}1$ induced by the $H$-coalgebra structure of $ML$
from (\ref{ML=H-coalgebra}).
\end{lemma}

\begin{proof}
Inductively. For $n=0$, there is nothing to show as $1$ is the terminal
object in $Set$. Assume $\alpha _{n}=a_{n} \circ Mp_{n}$, then in the following
diagram%
\begin{equation*}
\xymatrix{ ML \ar[r]^{Mp_{n+1}} \ar[d]_{M\xi} & MH^{n+1}1 \ar[r]^{{\lambda}_{H^n1}} & HMH^n1 \ar[d]^{Ha_n} \\
MHL \ar[ur]_{MHp_n} \ar[r]_{{\lambda}_L} & HML \ar[ur]^{HMp_n} \ar[r]_{H\alpha_n} & H^{n+1}1
}
\end{equation*}%
the triangle on the left commutes by (\ref{projections from the limit}), the
middle diagram commutes by naturality of $\lambda $ and the triangle on the right by
applying $H$ to the inductive hypothesis. It follows that $\alpha
_{n+1}=a_{n+1}\circ Mp_{n+1}$.
\end{proof}

In consequence, the unique coalgebra map $\gamma :ML\longrightarrow L$
constructed in (\ref{ML=H-coalgebra}) is also the anamorphism $\alpha
_{ML}:ML\longrightarrow L$ for the coalgebra $ML$.

\begin{lemma}
The projections $p_{n}:L\longrightarrow H^{n}1$ are $\mathbf{M}$-algebra
morphisms, with (\ref{ML=H-coalgebra}) and (\ref{Hn1=M-algebra}) giving the algebra
structures of $L$, respectively $H^{n}1$.
\end{lemma}

\begin{proof}
Again by induction. The first step is trivial. Assume that $p_{n}$ is an
algebra map: $\pi _{n} \circ \gamma =a_{n} \circ Mp_{n}$; then we have the following
diagram%
\begin{equation*}
\xymatrix{ ML \ar [dddd]_{\gamma} \ar [rr]^{Mp_{n+1}} \ar[dr]_{M\xi} \ar @{} [dddr] |{(2)}\ar @{} [drr] |{(1)}&  
& MH^{n+1}1 \ar [dd]^{\lambda_{H^n1}} \\
& MHL \ar[ur]_{MHp_n} \ar[d]^{\lambda_L} \ar@{} [dr] |{(4)}& \\
& HML \ar[d]^{H \gamma} \ar[r]^{HMp_n} & HMH^n1 \ar [dd]^{Ha_n}\\
& HL \ar @{} [ur] |{(5)}\ar[dr]^{Hp_n} & \\
L \ar @{} [urr] |{(3)} \ar[ur]^{\xi} \ar [rr]_{p_{n+1}} & & H^n1
}
\end{equation*}%
where:
(1) commutes by applying $M$ to (\ref{projections from the limit}); (2) commutes by (\ref{ML=H-coalgebra}); (3) commutes by (\ref{projections from the limit}); (4) commutes by the naturality of $\lambda $ and (5) commutes by applying $H$ to the inductive hypothesis.
\end{proof}

Resuming all above, we have the following diagram of $\mathbf{M}$-algebras
and $\mathbf{M}$-algebra morphisms, in which the lower sequence is limiting:
\begin{equation*}
\xymatrix{ M1 \ar[d]_{a_0} & MH1 \ar[l]_{Mt} \ar[d]_{a_1} & \ldots \ar[l] & MH^n1 \ar[l] \ar[d]_{a_n} & MH^{n+1}1 \ar[l]_{MH^nt} \ar[d]_{a_{n+1}} & \ldots \ar[l] & 
ML \ar@/_1.7pc/[lll]_{Mp_n} \ar@{.>}[d]^{\gamma} \\
1 & H1 \ar[l]_t & \ldots \ar[l] & H^n1 \ar[l] & H^{n+1}1 \ar[l]_{H^nt} & \ldots \ar[l] & L \ar@/^1.7pc/[lll]^{p_n}
}
\end{equation*}

\subsection{Topology on the final coalgebra}

From now on, we shall assume that $H$ is an $\omega^{op}$-continuous
endofunctor which admits a lifting to $Alg({\mathbf{M}})$. Remember that on all $H^{n}1$ we have consi\-dered the discrete topology.
Endow also all $MH^{n}1$ with the discrete topology (intuitively, this
corresponds to the fact that operations on algebras with discrete
topology are automatically continuous) and $ML$ with the initial
topology coming from the cone $Mp_{n}:ML\longrightarrow MH^{n}1$
(which is the same as the initial topology from the cone
$ML\overset{Mp_{n}}{\longrightarrow }%
MH^{n}1\overset{a_{n}}{\longrightarrow }H^{n}1$, as $a_{n}$ are
continuous maps between discrete spaces).

\begin{proposition}
  \label{topology on final coalgebra}Under the above assumptions, the
  final $H$%
  -coalgebra inherits a structure of a topological
  $\mathbf{M}$-algebra\footnote{Usually the notion of a topological algebra refers to an algebra for some finitary, algebraic theory whose underlying set is equipped with some topology, such that the algebra operations are continuous (\cite{stone}). As Eilenberg-Moore algebras for a $Set$-monad are the same as algebras for (not necessarily) finitary algebraic theories (\cite{ahs}), we find that the term topological algebra characterizes the best the present situation.}, i.e. $L$ has a $\mathbf{M}$-algebra structure $%
  \gamma :ML\longrightarrow L$ such that $\gamma$ is continuous with respect to the topologies on $L$ and $ML$.
\end{proposition}

\begin{proof}
  By definition of the initial topology, $\gamma $ is
  continuous if and only if all compositions $\gamma\circ p_{n}$ are
  continuous. But $\gamma\circ p_{n}=a_{n}\circ Mp_{n}$, $a_{n}$ are continuous
  as maps between discrete sets and $Mp_{n}$ are continuous by the
  initial topology on $ML$.
\end{proof}
 
Notice that this result relies heavily on the construction of the final
coalgebra as the limit of the sequence (\ref{terminal sequence}). Without
it, we can not obtain \linebreak Proposition \ref{topology on final coalgebra} just by assuming the existence of the final $H$-coalgebra and of the lifting to $Alg({\mathbf{M}})$, as there is no obvious choice
for the topology on $ML$. Also Proposition \ref{topology on final coalgebra} can be interpreted by saying that all
operations on $L$ are continuous (as they are obtained as limits of
operations on discrete algebras).

\begin{remark}
Instead of an $\omega^{op}$-continuous endofunctor, we could use a finitary
one. It is known \cite{worrell} that the final coalgebra exists, but the
previous limit is not enough. From this, a supplementary construction gives the final coalgebra. Obviously, the final coalgebra has an $\mathbf{M}$-algebra structure as in (\ref{ML=H-coalgebra}). Following Worrell's construction \cite{worrell}, the terminal sequence would still induce a topology on $L$, and the easiest way would be to take on $ML$ the initial topology with respect to $\gamma$, but this is not the same as the construction pursued here (the topology on $ML$ comes from the terminal sequence). 
\end{remark}

\subsection{Initial \~H-algebra and final \~H-coalgebra in Alg(\textbf{M})}

If $\widetilde{H}$ preserves colimits of $\omega $%
-sequences, then the initial $\widetilde{H}$-algebra is easy to build, using a dual procedure to the one in (\ref{terminal sequence}): recall that
$Alg({\mathbf{M}})$ has an initial object, namely the free algebra on
the empty set, \linebreak $F^{\mathbf{M}}0=(M0,M^20{\overset{m_0}{\longrightarrow}}M0)$. In order to simplify the notation, we shall
identify all algebras $\widetilde{H}^{n}F^{\mathbf{M}}0$ with their
underlying sets $%
H^{n}M0$. Then it is well-known that the initial
$\widetilde{H}$-algebra is the colimit in $Alg({\mathbf{M}})$ of the
chain
\begin{equation}{\label{initial sequence}}
M0\overset{!}{\longrightarrow }HM0\overset{H!}{\longrightarrow }%
...\longrightarrow H^{n}M0\overset{H^{n}!}{\longrightarrow }...
\end{equation}
where $!:M0\longrightarrow HM0$ is the unique algebra map. Denote by \linebreak
$i_{n}:H^{n}M0\longrightarrow I$ the colimiting cocone. We do not
detail anymore this construction as we did for coalgebras as it will
not be used in the sequel. However, we shall need the following (which
requires only the existence in $Alg(\mathbf{M})$ of the limit of the terminal sequence (\ref{terminal sequence}), respectively of the colimit of the initial sequence (\ref{initial sequence})): there is a unique \linebreak $\mathbf{M}$-algebra morphism 
$f:I\longrightarrow L$ such that
\begin{equation}
\xymatrix{ H^{n}M0 \ar[r]^-{i_{n}} \ar[d]_{H^{n}s} &I \ar[d]^{f} \\ 
H^{n}1 & L \ar[l]_-{p_{n}} }
\label{map colim-->lim}
\end{equation}%
commutes for all $n$ (see for example \cite{adamek ideal}, Lemma II.5 for a
proof), where \linebreak $s:M0 \longrightarrow 1$ is the unique algebra map from the initial to the final object in $Alg(\mathbf{M})$. If $M0$ not empty, then $I$ will also be not empty, as it comes with a cocone of algebra maps with not empty domains.

We shall generalize in this section the result of Barr (\cite{barr}) from $%
Set$ to $Alg({\mathbf{M}})$, for the special case of $Alg({\mathbf{M}})$%
-endofunctors arising as liftings of $Set$-endofunctors. The proofs use
similar ideas to the ones in \cite{barr} and \cite{adamek ideal}. 

We shall assume that there is an algebra map
\begin{equation}\label{zero object}
j:1 \longrightarrow M0  
\end{equation}
\noindent
As $M0$ is initial, $j \circ s=Id$. By finality of $1$ in $Alg(\mathbf{M})$, $s \circ j=Id$, hence we may identify $M0$ and $1$ as the zero object in the category of algebras.

\begin{remark}

There is a large class of monads satisfying this condition: the list monad
(and the commutative monoid-group-semi-ring monad), the (finite) power-set monad, the maybe monad, the $\Bbbk$-modules monad for a semi-ring $\Bbbk $. For all these, the free algebra with empty generators is built on the singleton set. But there are also monads for which the carrier of the free algebra on the empty set has more than one element, as the exception monad or the families monad, or it is empty, as is the case for the monad $MX=X \times \mathfrak{M}$, for $\mathfrak{M}$ a monoid. It is still under work whether the results of the present paper hold under this weakened assumption.
\end{remark}

We have $!:1=M0\longrightarrow HM0=H1$ and $t \circ !=Id$ in $Alg({\mathbf{M}})$. Hence
in the final sequence (\ref{terminal sequence}) all morphisms are split
algebra maps, the colimit is the initial $\widetilde{H}$-algebra and the
limit is the final $H$ (and $\widetilde{H}$)-coalgebra:

\begin{equation}
1\overset{t}{\underset{!}{\leftrightarrows }}H1\leftrightarrows
...\leftrightarrows H^{n}1\overset{H^{n}t}{\underset{H^{n}!}{%
\leftrightarrows }}H^{n+1}1\leftrightarrows ...
\label{split terminal sequence}
\end{equation}

\begin{theorem}
  \label{th completion initial final}Let $H$ a $Set$-endofunctor
  $\omega ^{op}$%
  -continuous, $\mathbf{M}$ a monad on $Set$ such
  that:

\begin{enumerate}
\item $H$ admits a lifting $\widetilde{H}$ to $Alg({\mathbf{M}})$ which is $%
\omega $-cocontinuous;

\item $M0=1$ in $Alg({\mathbf{M}})$;
\end{enumerate}
\noindent
then the final $H$-coalgebra is the completion of the initial $\widetilde{H}$-
algebra under a suitable (ultra)metric.
\end{theorem}

\begin{proof}
Consider the following diagram (in $Alg({\mathbf{M}})$), where all algebras
involved have structure maps defined via the distributive law $\lambda $.%
\begin{equation*}
   \xymatrix{
   1 \ar@<0.5ex>[r]^{!} & H1 \ar@<0.5ex>[l]^t \ar@<0.5ex>[r] & \ldots \ar@<0.5ex>[l] \ar@<0.5ex>[r] & H^n1 \ar[dl]^{i_n} \ar@<0.5ex>[l] \ar@<0.5ex>[r]^{H^{!}} & \ldots \ar@<0.5ex>[l]^{H^t} \\
   & & I \ar[rr]_f & & L \ar[ul]^{p_n} 
   }
\end{equation*}
\noindent
Put on $I$ the smallest topology such that $f$ is continuous, where $L$ has the structure of a topological algebra from Proposition \ref{topology on
final coalgebra}. This coincides with the initial topology given by the cone 
$I\overset{f}{\longrightarrow }L\overset{p_{n}}{\longrightarrow }H^{n}1$.
Moreover, $I$ becomes a topological algebra and all $i_{n}$ are continuous algebra
maps, if on $MI$ we take the topology induced by the map $%
Mf:MI\longrightarrow ML$. In particular, $Mf$ is continuous. Denote by $MI%
\overset{\zeta }{\longrightarrow }I$ the algebra structure map of $I$. Then $%
f \circ \zeta =\gamma \circ Mf$ (remember that $f$ is an algebra map). As $L$ is a
topological $\mathbf{M}$-algebra, it follows that $f \circ \zeta $ is continuous,
hence $\zeta $ is continuous.
About $i_{n}$: these are by construction algebra maps (as the components of the colimiting cocone in $Alg(\mathbf{M})$) and also continuous, as $H^{n}1$ are discrete.
The only remaining thing we need to prove is the density of $I$ (more precisely, of $Imf$) in $L$. We start by applying Barr's argument to show that $L$ is
complete under this ultrametric. First, use that limits in $Alg({\mathbf{M}})$
are computed as in $Set$ to conclude that $L$ is Cauchy complete: take a
Cauchy sequence $x^{(n)}$ in $L$ with respect to the initial topology (ultrametric) and
assume $d(x^{(n)},x^{(m)})<2^{-\min (m,n)}$ for all $m$, $n$. This implies $%
p_{n}\circ f(x^{(n)})=p_{n}\circ f(x^{(m)})$ for all $n<m$. Thus $y=(p_{n}\circ f(x^{(n)}))_{n%
\geq 0}$ defines an element of $L$ such that $\lim x^{(n)}=y$. Next, a similar construction to the one in \cite{adamek
final} will show us that the image of $I$ under the algebra morphism $f$ is
dense in $L$. For this purpose, consider the additional $\mathbf{M}$-algebra
sequence of morphisms $(h_{n})_{n\geq 0}$, given by%
\begin{equation*}
h_{n}:L\overset{p_{n}}{\longrightarrow }H^{n}1=H^{n}M0\overset{H^{n}!}{%
\longrightarrow }H^{n+1}M0\overset{i_{n+1}}{\longrightarrow }I\overset{f}{%
\longrightarrow }L
\end{equation*}

We have $p_{n+1} \circ h_{n}=H^{n}!\circ p_{n}$. Consider now an element $x\in L$. Then
by construction $(y^{(n)}=h_{n}(x))_{n\geq 0}$ form a sequence of elements
lying in the image of $f$ and we shall see that this sequence is convergent
to $x$. Indeed, from $p_{n+1}(y^{(n)})=H^{n}! \circ p_{n}(x)$ it follows that 
\begin{equation*}
p_{n}(y^{(n)})=H^{n} \circ t \circ p_{n+1}(y^{(n)})=H^{n} \circ t \circ H^{n}!\circ p_{n}(x)=p_{n}(x)
\end{equation*}%
the $n$-th projection of the $n$-th term of the sequence $(y^{n)})_{n\geq 0}$
coinciding with the $n$-th projection of the element $x$; hence $%
d(y^{(n)},x)<2^{-n}$ which implies $\lim y^{(n)}=x$
in $L$. Therefore the image of $I$ through the canonical colimit$\longrightarrow$limit arrow is dense in $L$.
\end{proof}

\begin{remark}\label{a}
\begin{enumerate}
\item If we consider on the initial algebra $I$ the final topo\-logy coming
from the $\omega $-chain, this is exactly the discrete topology (and \linebreak
metric), since all $H^{n}1$ are discrete, hence $I$ would be Cauchy complete
and $f:I\rightarrow L$ automatically continuous. No interesting information
between $I$ and $L$ can be obtained in this situation.
\item From (\ref{map colim-->lim}) and (\ref{zero object}) we have $p_{n}\circ f\circ i_{n}=Id$, hence $f\circ i_{n}$ is a monomorphism. But all morphism in the above sequence are split algebra maps by (\ref{split terminal sequence}), hence all $H^n{!}$ are mono's. Recall now from \cite{adamek ideal} that in any locally finitely presentable category, 
\begin{itemize}
\item the cocone to the colimit of an $\omega $-chain formed by monomorphisms is a monomorphism,
\end{itemize}
and 
\begin{itemize}
\item for every cocone to the chain formed by monomorphisms, the unique
map from the colimit is again a monomorphism.
\end{itemize}
If we assume $M$ finitary, the Eilenberg-Moore category of algebras would be locally finitely presentable. Hence the algebra map $f$ would be mono. But remember that any $Set$-monad is regular (\cite{ahs}). It follows that we can identify $I$ with a subalgebra of $L$. The algebra isomorphism $g:I\simeq Imf$ would also be a homeomorphism, if we take on $Imf$ the induced topology from $L$. 
\item \label{remark} The $\omega$-cocontinuity of $\widetilde{H}$ is automatically
  satisfied if we assume $M,H$ to be finitary. For, the monad being finitary,
  the forgetful functor $U^{\mathbf{M}}$ would preserve and reflect
  sifted colimits. But $U^{\mathbf{M}}{\widetilde{H}}=HU^{\mathbf{M}}$,
  hence $\widetilde{H}$ commutes with sifted colimits, in particular with
  colimits of $\omega$-chains.
\end{enumerate}
\end{remark}

\begin{example}\label{semi-ring} The functor $HX=\Bbbk \times X^{A}$ is built from products, hence is \linebreak $\omega$-continuous. The $H$-coalgebras are known as Moore automata. Such a functor always admits at least one lifting to $Alg(\mathbf{M})$ for any monad $\mathbf{M}$, provided $\Bbbk$ carries an algebra structure. The lifted functor is given by the same formula as $H$, where this time the product and the power are computed in the category of algebras. 

In particular, consider $A$ a finite set, $\Bbbk $ a (not necessarily commutative) semi-ring and $\mathbf{M}$ the monad that it induces (as in \cite{maclane}, Section VI.4, Ex. 2, where the ring $R$ is replaced by the semi-ring $\Bbbk $);
then $Alg({\mathbf{M}})$ is the category of $\Bbbk $-modules and $M0$ is the
zero module. The final $H$-coalgebra is $\Bbbk ^{A^{\ast }}$, the set of all functions $A^{\ast}\longrightarrow \Bbbk$, also known as the formal power series in non-commuting $A$ variables, while the initial $\widetilde{H}$-algebra is the direct sum of $A^{\ast }$ copies of $\Bbbk $ (the polynomial algebra in same variables) (recall that in this case, finite products and coproducts coincide in $Alg(\mathbf{M})$). The approximants of order $n$ in the corresponding $\omega $-sequence are $H^{n}1=\Bbbk ^{1+A+...+A^{n}}$, the polynomials in (non-commuting) $A$-variables of degree at most $n$. We shall detail this for the easiest case, where $A$ is the singleton $\{t\}$; the distance between two elements of the final coalgebra $\Bbbk [[t]]$, i.e. between two power series $f(t),g(t)$ in variable $t$, is given precisely by $2^{-ord(f(t)-g(t))}$, where $ord(f(t)-g(t))$ is the order of the difference $f(t)-g(t)$ (the smallest power of $t$ which occurs with a nonzero coefficient in the difference). Take a Cauchy sequence of \linebreak polynomials $f_n(t)=a_{0}^{n} +a_{1}^{n}t+ \ldots$, where only finitely many $a_{j}^{n}$ are nonzero, for each $n,j \in \mathbb{N}$. For every $r\geq 0$, there exists an $n_r$ such that for every $n\geq n_r$, we have $ord(f_n(t)-f_{n_r} (t))=r$; this implies $a_{j}^{n}=a_{j}^{n_r}$ for all $j\leq r$ and $n \geq n_r$. Let $f(t)=a_{0}^{n_0}+a_{1}^{n_1}t+ \ldots$. One immediately verifies that the power series $f(t)$ is the limit of of the sequence $(f_n(t))_{n\geq 0}$. Hence the final coalgebra $\Bbbk [[t]]$ is indeed the completion of the initial $\widetilde{H}$-algebra $\Bbbk [t]$. 
\end{example}

\section{Application: \textbf{M}-commuting pairs of endofunctors}

Consider an endofunctor $H$ and a monad $\mathbf{M}$, both on $Set$. There are two ways of relating the endofunctor to the monad by a natural transformation, as follows:
\begin{itemize}
\item 
$\lambda : MH \longrightarrow HM$ satisfying (\ref{distributive law endofunctor monad}), which is the same as an algebra lift $\widetilde{H}:Alg({\mathbf{M}}) \longrightarrow Alg({\mathbf{M}})$, $U^{\mathbf{M}}{\widetilde{H}}=HU^{\mathbf{M}}$;
\end{itemize}
or
\begin{itemize}
\item
$\varsigma :HM \longrightarrow MH$ satisfying \\
\begin{allowdisplaybreaks}
\begin{equation}\label{kleisli lift}
	\xymatrix{ H \ar[r]^{Hu} \ar[dr]_{u_H} & HM \ar[d]^{\varsigma} \\
	 & MH}
\qquad
	\xymatrix {HM^{2} \ar[r]^{{\varsigma}_M} \ar[d]_{Hm} & MHM \ar[r]^{M\varsigma} & 
	M^2H \ar[d]^{m_H} \\
	HM \ar[rr]^{\varsigma} & & MH}
\end{equation} 
\end{allowdisplaybreaks}
\noindent
It is well known that this is equivalent to the existence of a Kleisli lift, i.e. an endofunctor $\hat{H}:Kl({\mathbf{M}}) \longrightarrow Kl({\mathbf{M}})$ such that $\hat{H}F_{\mathbf{M}}=F_{\mathbf{M}}H$, where \linebreak $F_{\mathbf{M}}:Set \longrightarrow Kl({\mathbf{M}})$ is the canonical functor to the Kleisli category of the monad. In this case, we can perform the following additional construction: denote by 
$\mathcal{I} : Kl({\mathbf{M}}) \longrightarrow Alg({\mathbf{M}})$ the comparison functor. Take the $Alg({\mathbf{M}})$-endofunctor given by the left Kan extension along $\mathcal{I}$ (which exists since every algebra in $Alg({\mathbf{M}})$ arises as a coequaliser of
free algebras in a canonical way):
\begin{equation}\label{kleisli} 
\bar{H}=Lan_{\mathcal{I}}({{\mathcal{I}}{\hat{H}}})
\end{equation}
\noindent
As the Kleisli category $Kl(\mathbf{M})$ is isomorphic to a full subcategory of $Alg(\mathbf{M})$, this would yield a natural isomorphism $\mathcal{I} \hat{H} \cong \bar{H} \mathcal{I}$. Composing this with the functor $F_{\mathbf{M}}$, we obtain  
$\bar{H}F^{\mathbf{M}} \cong F^{\mathbf{M}}H$, as in the diagram below:
\begin{equation}
\xymatrix{
Alg(\mathbf{M}) \ar@{.>}[r]^{\bar{H}} & Alg(\mathbf{M}) \\
Kl(\mathbf{M}) \ar[u]_{\mathcal{I}} \ar[r]^{\hat{H}} & Kl(\mathbf{M}) \ar[u]^{\mathcal{I}} \\
Set \ar@/^2pc/[uu]^{F^{\mathbf{M}}} \ar[u]_{F_{\mathbf{M}}} \ar[r]^H & Set \ar@/_2pc/[uu]_{F^{\mathbf{M}}} \ar[u]^{F_{\mathbf{M}}}
}
\end{equation}
We shall call $\bar{H}$ an extension of $H$ to algebras. 
\end{itemize}

With the above notations, consider now two $Set$-functors $T$, $H$ such that both an algebra lift of $H$ and a Kleisli lift of $T$ exist and $\widetilde{H} \cong \bar{T}$. Then we have 
\begin{eqnarray*}
MT &=&U^{\mathbf{M}}F^{\mathbf{M}}T\cong U^{\mathbf{M}}{\bar{T}}F^{\mathbf{M}}
\\
& \cong &U^{\mathbf{M}}{\widetilde{H}}F^{\mathbf{M}}=HU^{\mathbf{M}}F^{\mathbf{M}}=HM
\end{eqnarray*}
\noindent
i.e. $M$ acts like a switch (up to isomorphism) between the endofunctors $T$ and $H$. 
\begin{definition}
Let $\mathbf{M}=(M,m,u)$ be a monad on $Set$. A pair of $Set$-endofunctors $(T,H)$ such that $HM \cong MT$ is called an $\mathbf{M}$-commuting pair. 
\end{definition}

\begin{example}
One can easily obtain commuting pairs in the following situations:
\begin{itemize}
\item Take $T=H=Id$ or $T=H=M$ and $\mathbf{M}=(M,m,u)$ any monad;
\item Consider $T=H=A+(-)$, $\mathbf{M}=B+(-)$. Then commmutativity of the coproduct ensures the commuting pair; similarly for products: $T=H=A \times (-)$, $\mathbf{M}=B \times (-)$, where this time $B$ is a monoid (this works more generally, in any monoidal category).
\end{itemize}

\end{example}
To the best of our knowledge, it seems that the notion of commuting pairs has not been considered previously, although the above examples show that it arises naturally in mathematics. We shall later see more (non-trivial) examples. But before that, we come back to the situation considered earlier, of the two endofunctors $T$ and $H$ such that $\widetilde{H} \cong \bar{T}$. This implies 
\begin{equation*}
\widetilde{H}F^{\mathbf{M}} \cong \bar{T}F^{\mathbf{M}} \cong F^{\mathbf{M}}T
\end{equation*}
\noindent 
which can be rephrased by saying that $HM \cong MT$ is an isomorphism of $\mathbf{M}$-algebras, where the algebra structure of $HMX$, for a set $X$, is induced by the distributivity law $\lambda:MH \longrightarrow HM$, i.e. the following diagram commutes:
\begin{equation}\label{algebra isom}
\xymatrix{
MHMX \ar[r]^{\cong} \ar[d]_{\lambda_{MX}} & M^2TX \ar[dd]^{m_{TX}} \\ HM^2X \ar[d]_{Hm_X} \\ HMX \ar[r]^{\cong} &
MTX
}
\end{equation}
\noindent where the lower horizontal arrow is $HM \cong MT$, while the upper arrow is obtained by applying $M$ to this.

Conversely, if $(T,H)$ is an $\mathbf{M}$-commuting pair, one may
wonder about their relation with the category of
$\mathbf{M}$-algebras. Suppose $H$ has an algebra lifting $\widetilde{H}$,
$T$ has a Kleisli lift (hence an extension $\bar{T}$) and $HM \cong MT$ such that (\ref{algebra isom}) holds; then from $HM \cong MT$ and 
\begin{eqnarray*}
& HM=HU^{\mathbf{M}}F^{\mathbf{M}}=U^{\mathbf{M}} \widetilde{H} F^{\mathbf{M}} \\
& MT=U^{\mathbf{M}}F^{\mathbf{M}}T \cong U^{\mathbf{M}} \bar{T} F^{\mathbf{M}}
\end{eqnarray*}
it follows that $U^{\mathbf{M}} \widetilde{H} F^{\mathbf{M}} \cong
U^{\mathbf{M}} \bar{T} F^{\mathbf{M}}$, that is, the images of
$\widetilde{H}$ and $\bar{T}$ on free algebras share (up to bijection)
the same underlying set. Taking into account that $HM \cong MT$ is an
isomorphism of $\mathbf{M}$-algebras (\ref{algebra isom}), we obtain
that $\widetilde{H} \cong \bar{T}$ on free algebras. Assume now that
$M$, $T$ and $H$ are finitary. Then, by construction, $\bar{T}$ is
determined by its action on finitely generated free algebras, and so
is $\widetilde H$ (because it preserves sifted colimits by Remark
\ref{a}(iii)). It follows that $\widetilde H \cong \bar T$.

We have obtained thus

\begin{proposition}\label{commuting}
  Let $H$, $T$ two endofunctors on $Set$ and $\mathbf{M}$ a monad on
  $Set$. Assume that $H$ has an algebra lift $\widetilde{H}$ and $T$ has a
  Kleisli lift with respect to the monad $\mathbf{M}$. Denote by
  $\bar{T}$ the corresponding left Kan extension, as in
  (\ref{kleisli}). Then:
\begin{enumerate}
\item If $\widetilde{H} \cong \bar{T}$, then $(T,H)$ form an $\mathbf{M}$-commuting pair and $HM \cong MT$ is an algebra isomorphism.
\item Conversely, if $M,H,T$ are finitary and $MT \cong HM$ as algebras, then $\widetilde{H} \cong \bar{T}$.
\end{enumerate}
\end{proposition}

\begin{example}
  \label{exemplu}Take $TX=1+A\times X$, with $A$ finite set and
  $\mathbf{M}$ any $Set$-monad. Then a Kleisli lifting of $T$ exists, namely for each map $X \overset{f}{\longrightarrow} MY$, take $TX \overset{f}{\longrightarrow}MTY$ to be the composite 
\begin{eqnarray*}
&&TX=1+A\times X \overset{1+A \times f}{\longrightarrow}1+A \times MY\longrightarrow \\
&&1+M(A \times Y)\longrightarrow M1+M(A \times Y)\longrightarrow M(1+A \times Y)
\end{eqnarray*}
where the map $1+A \times MY\longrightarrow 1+M(A \times Y)$ is obtained from the canonical strength of the monad, while $1+M(A\times Y)\longrightarrow M1+M(A \times Y)$ uses the unit of the monad and $M1+M(A \times Y)\longrightarrow M(1+A \times Y)$ comes from the coproduct property. Also, it is easy to see that the extension of $T$ to $\mathbf{%
    M}$-algebras is $\bar{T}X=F^{\mathbf{M}}1+A\cdot X$, for each algebra $%
  X$, where this time the coproduct (respectively the copower) is
  computed in $Alg({\mathbf{M}})$. If the category of
  $\mathbf{M}$-algebras has finite biproducts (as in the case of the monad induced by a semi-ring, see Example \ref{semi-ring}), then $\bar{T}$ is the lifting to $Alg(\mathbf{M})$ of the $Set$-endofunctor $HX=M1\times X^{A}$. Hence $(T,H)$ form a commuting pair.
\end{example}
The motivation for studying commuting pairs appears clearly if we
combine the previous proposition with our main result from Theorem
\ref{th completion initial final}, obtaining the following:
\begin{corollary}\label{corol}
  Assume the assumptions of Proposition \ref{commuting}(ii) hold. If
  $H$ is $\omega ^{op}$-continuous and $M0=1$ as
  $\mathbf{M}$-algebras, then the final $H$-coalgebra is the
  completion of the free $\mathbf{M}$-algebra built on the initial
  $T$-algebra under a suitable metric.
\end{corollary}

\begin{proof}
  Follows from Theorem \ref{th completion initial final}, by noticing
  that the $M$-image of the initial $T$-algebra (which exists as $T$
  is finitary, hence $\omega $-cocontinuous) is the initial
  $\bar{T}$-algebra (by construction, $\bar{T}$ is finitary, so
  $\omega $-cocontinuous), while $ H$ and $\widetilde{H}$ share same final
  coalgebra.
\end{proof}

\begin{example}
We come back to Example \ref{exemplu} and take the monad induced by a semi-ring $\Bbbk$, as in Example \ref{semi-ring}. Then the initial $T$%
-algebra is $A^{\ast }$, the monoid of all finite words (including the empty one)
built on the alphabet $A$. The free $%
\mathbf{M}$-algebra is the direct sum of $A^{\ast }$
copies of $\Bbbk $, that is, the polynomial $\Bbbk $-algebra in non-commuting
$A$-variables $\Bbbk \lbrack A]$ (in the category of $\Bbbk $-semimodules),
while the final $H$-coalgebra is $\Bbbk ^{A^{\ast }}$, the
non-commutative power series $\Bbbk $-algebra.
\end{example}
The situation described until now in this section can be presented as
follows: If two endofunctors $T$ and $H$ are given, one may search for
the appropriate monad such that $(T,H)$ form a commuting pair. As
there is a special bond between algebras of $T$ and coalgebras of $H$,
it is not clear whether the general case of any two (finitary)
$Set$-endofunctors would have a solution. But there is another
possible approach: Start only with one endofunctor and additionally
with a (finitary) monad; find then a distributive law inducing a
Kleisli (or algebra) lift. Once this is accomplished, one should built
a second endofunctor on $Set$ (assuming this is possible) in order to
obtain a commuting pair, using the functor obtained on
$Alg(\mathbf{M})$.

For lifting to the Kleisli category, there is the following suitable
situation: for all commutative monads $\mathbf{M}$ and all analytic
functors $T$, a distributive law $TM\longrightarrow MT$ can always be
constructed (\cite{milius}). The commutativity of $\mathbf{M}$ ensures
also the existence of a tensor product $\otimes$ on $Alg(\mathbf{M})$,
such that the free functor $F^{\mathbf{M}}:(Set,\times)
\longrightarrow (Alg(\mathbf{M}), \otimes)$ is strong monoidal
(\cite{weakening}). If $T$ is a polynomial functor
$TX=\coprod\limits_{n\geq0} A_n \times X^n$, an obvious choice of
Kleisli lift would give (the extension)
${\bar{T}}X={\coprod\limits_{n\geq0}}F^{\mathbf{M}}A_n \otimes
X^{\otimes n}$, where this time $X \in Alg(\mathbf{M})$. Now recall
that both the coproduct and the tensor product on $Alg(\mathbf{M})$
are obtained as reflexive coequalizers, hence if we assume the monad
not only commutative but also finitary (as all results in this section
rely on the finitariness of $M$), it follows that the forgetful
functor would transform the coproduct, respectively the tensor product
of any two algebras $(X,x)$, $(Y,y)$ into a reflexive coequalizer
computed this time in $Set$. In particular, for the polynomial functor
$T$, a corresponding commuting pair $(T,H)$ exists and can be
constructed by the above argument. Moreover such $H$ is finitary by
construction. If $H$ is also ${\omega}^{op}$-continuous and $M0=1$ as
algebras, then by Corollary \ref{corol} the final $H$-coalgebra should
be realized as a completion of the (image) of the free algebra built
on the initial $T$-algebra (which is well known to be the set of
finite trees with branching and labeling given by the signature of
$T$).

However, lifting functors to the Eilenberg-Moore category seems to be
more problematic, even for the simplest case of polynomial functors, as follows:
\begin{itemize}
\item if $H$ is a constant functor, then the image of $H$ (the set) must be
the carrier of an $\mathbf{M}$-algebra $(A,a)$; if this is the case, one may form a commuting pair if and only if $A$ is a free algebra. Then $T$ is also a constant functor; in particular, Corollary \ref{corol} is trivially true.  

\item if $HX=A\times X$, and $A$ is the carrier of an algebra, a lift is easily seen to exist, as the forgetful functor $U^{\mathbf{M}}$ preserve products. Conversely, if $\widetilde{H}$ is a lifting of $H$, then there is an algebra structure on $A$, namely ${\widetilde{H}}1$. If the category $Alg(\mathbf{M})$ has finite biproducts (for example if $\mathbf{M}$ is the monad induced by a semi-ring $\Bbbk$) and $A$ is the carrier of a free algebra with set generators $B$, then there is a commuting pair $(T,H)$ with $TX=B+X$. The final $H$-coalgebra is the set of all streams on $A$, while the $\mathbf{M}$-algebra on the initial $T$-algebra is the $\omega$-copower of $MB\cong A$. 
\item if $HX=X^{n}$, a finite power functor, then the lifting exists as the
forgetful functor $U^{\mathbf{M}}$ preserves limits; the existence of finite biproducts in $Alg(\mathbf{M})$ is again the most convenient way of finding the correspondent functor as a copower $TX=n\cdot X$. But in this case no relevant answers are obtained in the initial-final (co)algebra relation, as these objects are trivial (empty initial $T$-algebra, singleton final $H$-coalgebra).

\item if $HX=A+X$ or $HX=X+X$, there is no obvious distributive law $\lambda :MH\longrightarrow HM$, unless the monad itself is obtained as a sum (like the maybe monad $MX=1+X$). 
\end{itemize}

\section{Conclusions}
The general picture behind Barr's theorem is conceptually simpler: if one starts with an arbitrary category $\mathcal{C}$ (with initial object, final object and \linebreak $\omega$-(co)limits) and a $\mathcal{C}$-functor, then the theorem roughly says that the $\omega$-limit of the terminal sequence is a completion of the ${\omega}^{op}$-colimit of the initial sequence. Of course an appropriate notion of completion is required; it could be of topological nature (as in \cite{barr}), or about ordered structures (\cite{adamek ideal}). In the present paper we have emphasized the topological aspect (Cauchy completion) for base category $\mathcal{C}$ with algebraic structure, namely the Eilenberg-Moore category of a $Set$-monad. The endofunctors considered were obtained as liftings from $Set$, as one of our motivations came from the following question: given a continuous $Set$-functor $H$ with $H0=0$, what can be said about the final $H$-coalgebra? If the functor is not necessary continuous (for example the finite powerset functor), then the final sequence has to be extended beyond $\omega$ steps. What happens with the completion procedure on $Alg(\mathbf{M})$ in such cases? We believe that an answer to this question is worth considering in the future. 

The second part of the paper introduces the notion of a commuting pair of endofunctors with respect to a monad. This seems to be new, however a detailed analysis and more examples are needed in order to better understand this structure (like the connection between bisimulations and traces exhibited in \cite{silva}). We plan to do this in a further paper.

\end{document}